\documentclass[reqno,12pt]{amsart}
\usepackage{amssymb,amsmath,amsthm,hhline,graphicx}

\newtheorem{theorem}{Theorem}

\newtheorem*{lemma}{Lemma}

\theoremstyle{definition}

\newcommand{\PSL}{\mathrm{PSL}}
\newcommand{\SL}{\mathrm{SL}}
\newcommand{\PGL}{\mathrm{PGL}}

\begin{document}

\title{Galois's lost proof on $\PSL(2,p)$}
\author{P\'eter P.~P\'alfy}
\address{Alfr\'ed R\'enyi Institute of Mathematics \\
Re\'altanoda utca 13--15, 1053 Budapest, Hungary}
\email{ppp@renyi.hu}
\dedicatory{Dedicated to the memory of Peter M.\ Neumann (1940--2020)}
\thanks{This research was supported by the National Research, Development and
Innovation Office of Hungary (NKFIH), grant no.~153681}
\keywords{\'Evariste Galois, projective special linear group, imprimitive permutation group}

\begin{abstract}
In his famous ``testamentary letter'' Galois formulated several results without proof.
Among them he claimed --- here we express it in modern terms ---
that the projective special linear group $\PSL(2,p)$ has a subgroup of index $p$ only for $p\le 11$.
Nowadays this result is seen as a corollary of the complete enumeration of all subgroups
of $\PSL(2,p)$ obtained by Gierster in 1881.
In this paper we attempt to give a proof that uses only methods available to Galois.
\end{abstract}
\maketitle

\section{Introduction}

The famous ``testamentary letter'' written by \'Evariste Galois on the eve of his fatal
duel contains a wealth of mathematical results. In the first part of the letter Galois
summarizes his theory on the solvability of polynomial equations what was the subject of
his ``First Memoir'' submitted to the French Academy on January 16, 1831.
(For a modern treatment of this work see \cite{Edwards}.)
After that almost a quarter of the letter
(namely, the bottom part of page 3, and the whole page 4) is devoted to the study of
the Galois group of the modular equation arising in the theory of elliptic functions.
In modern terms this group is the projective general linear group $\PGL(2,p)$
acting on the $p+1$ roots of the modular equation indexed by the elements of the
projective line over the $p$-element field ($p$ a prime number), i.e., $\{0,1,\dots,p-1,\infty\}$.
The action is given by fractional linear transformations $x\mapsto (ax+b)/(cx+d)$, where
$ad-bc\ne 0$. In particular, $\infty\mapsto a/c$ and $-d/c\mapsto\infty$.

Galois first notices that the group has a normal subgroup of index 2, namely, the
subgroup consisting of those fractional linear transformations, where $ad-bc$ is a
quadratic residue modulo $p$
(in other words, a nonzero square in the $p$-element field). Today this normal subgroup is
called the projective special linear group $\PSL(2,p)$. For $p=2$ and $p=3$ the degree of
the modular equation is 3 and 4, respectively, so the equation is solvable by radicals,
and the group is solvable. So we will assume that $p\ge 5$.

Galois has two substantial results concerning $\PSL(2,p)$. We will call them theorems,
although they are not designated so in the running text of the testamentary letter.

\begin{theorem}
Let $p\ge 5$ be a prime. The projective special linear group $\PSL(2,p)$ is a simple group.
\end{theorem}

\begin{theorem}
Let $p\ge 5$ be a prime. The projective special linear group $\PSL(2,p)$ has a subgroup of
index $p$ only when $p=5$, $7$, and $11$.
\end{theorem}

Concerning the first theorem he simply says ``But it is easy to see that it is not further
decomposable properly unless $p=2$ or $p=3$.'' (See \cite[p.~89]{Neumann}. The English translations
of Galois's texts are taken from this edition of Galois's works by Peter M.\ Neumann.)
We have no clue what kind of easy proof he had in mind.

I have always been admiring Galois for the second theorem. Because of the three exceptional
cases it is absolutely non-trivial. Until recently, I have been happy to see this result
as a consequence of the complete classification of subgroups of $\PSL(2,q)$ (where $q$ is
a prime power), see, e.g., \cite[Kapitel II, \S 8]{Huppert}.
But this is based on much later developments in group theory. So how could
Galois prove this result? Again, he does not say anything about the proof.
After outlining the examples for $p=5$, $7$, $11$ he states boldly:
``In all rigour this reduction is not possible in the higher cases.''
(See \cite[p.~91]{Neumann}.)
However, analysing his constructions for the small primes may provide hints for a
possible proof applying only tools already available to him. In this note we try to give
such a proof in Section~\ref{theproof}.

Before doing that we make some historical remarks in Section~\ref{historical}.

Finally, in Section~\ref{examples}, we relate Galois's examples for the cases
$p=5,7,11$ to the rotation group of the tetrahedron, the cube, and the
icosahedron, respectively.

\section{Historical remarks}\label{historical}

Except for some scattered calculations with fractional linear transformations, Galois's manu\-scripts do
not contain any explanation for his claims about $\PSL(2,p)$. The examples he gave for the cases
$p=5,7,11$ may give some hints how he looked at these subgroups. By describing a pairing of the
permuted elements (for details see the next section), he is obviously indicating that these subgroups
are imprimitive permutation groups. One finds the definition of imprimitivity
in his handwriting on a torn piece of paper among the fragments
in Dossier 15 at the Library of the French Academy (see \cite[p.~289]{Neumann}). Our suggested
reconstruction of Galois's proof is based on this concept.

One might wonder why Galois wrote ``In all rigour'' introducing the statement about the
nonexistence of subgroups of index $p$ for $p>11$. On the one hand, he tacitly corrects an
earlier error. Namely, in the short announcement of his results on the theory of equations
(see \cite[p.~53]{Neumann}) that was published in F\'erussac's Bulletin in April, 1830, he
erroneously stated that such a reduction (i.e., the existence of such a subgroup) is not
possible even for $p>5$. On the other hand, he is now absolutely certain that the result in
this new version is correct, as he writes on the last page of the testamentary letter: ``Often in my life
I have risked advancing propositions of which I was not certain. But all that I have written here
has been in my head for almost a year and it is not in my interest to make a mistake so that one
could suspect me of having announced theorems of which I did not have the complete proof''.
This sentence also suggests that he had not written down the proof of this result. The manuscript
of his ``Second Memoir'' (see \cite[Chapter V]{Neumann}) is unfinished, it does not reach the
point of our interest, leaving room for guesswork.

Enrico Betti \cite{Betti} was the first to publish a proof of Galois's result (Theorem 2 above)
in 1853.
However, as it is sometimes the case with Betti's proofs (cf.~\cite{Mammone}), his argument is
not complete. In fact, he gives a sufficient condition for the existence of a subgroup of
index $p$ in $\PSL(2,p)$, and checks that this is not satisfied by primes greater than 11.
Nevertheless, his intuition was right. The necessary condition we will derive in
Section~\ref{theproof} is quite similar to Betti's sufficient condition. Here is the English
translation (with my explanations in square brackets) of his theorem :

\begin{theorem}
(Betti \cite[p.~88]{Betti})
When the odd powers $<p-2$ of a primitive root $g$ of the prime number
$p=4n+3$, for which $\left( \frac{g-1}{p}\right)=1$, satisfy either one or the other
of the two congruences
\begin{gather}
g^2x^2-g(g+1)x+1\equiv 0 \pmod{p},\\
g^2x^2-(g+1)x+1\equiv 0 \pmod{p},
\end{gather}
then the group $M$ [$=\PSL(2,p)$] is the sum [i.e., union] of $p$ derived groups similar to joint letters
[i.e., cosets of a subgroup].
\end{theorem}

Charles Hermite \cite{Hermite}, \cite{HermiteOeuvres} was intrigued to find explicit equations of degree $p$ that
are arise as reductions of the modular equation of degree $p+1$ in the cases $p=5,7,11$.
(See also \cite{Goldstein}.)

The first correct proof of Theorem 2 can be found in Camille Jordan's
{\it Trait\'e des substitutions et des \'equations alg\'ebriques} published in 1870. He works in
the matrix group $\SL(2,p)$ and using linear algebraic methods he shows that this group
cannot have any subgroup of index $p$, whenever $p>11$. His general reasoning requires
$p>13$, so he has to add some ad hoc argument to deal with the case $p=13$. (See \cite[pp.~347--351]{Jordan}.)

In his doctoral thesis \cite{Gierster} Felix Klein's student, Joseph Gierster determined all subgroups
of $\PSL(2,p)$. Galois's theorem immediately follows from this result. More generally, subgroups of
the projective special linear groups $\PSL(2,q)$ over arbitrary finite fields of $q$ elements, have
been described by Dickson \cite{Dickson} (see also \cite{Huppert}).

\section{The proof}\label{theproof}

For preparation we first prove a simple property of the group of linear functions
over $\mathbb{F}_p$, the $p$-element field (where $p$ is a prime number).

\begin{lemma}
Any subgroup in the group of linear functions over the $p$-element field
whose order is not divisible by $p$ fixes an element.
\end{lemma}

\begin{proof}
Let $L$ be a subgroup in the group of linear functions $x\mapsto ax+b$, where
$x,a,b\in\mathbb{F}_p$, $a\ne 0$, and assume that $p$ does not divide $|L|$.
If $a=1$, $b\ne 0$, then $x\mapsto x+b$ has order $p$, so such a linear function
cannot belong to $L$. If $a\ne 1$, then $x\mapsto ax+b$ has a unique fixed point
$b/(1-a)$. We have to show that if $f_1(x)=a_1x+b_1$ with $a_1\ne 1$ and
$f_2(x)=a_2x+b_2$ with $a_2\ne 1$ both belong to $L$, then they have the same
fixed point, i.e., $b_1/(1-a_1)=b_2/(1-a_2)$. Indeed, we have
$f_1(f_2(x))=a_1(a_2x+b_2)+b_1=a_1a_2x+(a_1b_2+b_1)$, and similarly
$f_2(f_1(x))=a_1a_2x+(a_2b_1+b_2)$, hence $f_1f_2=gf_2f_1$, where
$g(x)=x+(a_1b_2+b_1)-(a_2b_1+b_2)$. Here $g=(f_1f_2)(f_2f_1)^{-1}\in L$, thus
by the observation above we must have $(a_1b_2+b_1)-(a_2b_1+b_2)=0$, that is
$b_1/(1-a_1)=b_2/(1-a_2)$, as it was required.
\end{proof}

Now let's turn to the proof of Galois's theorem. The constructions
in the exceptional cases will be given in the next section.

\begin{proof}[Proof of Theorem~2]
Let $G=\PSL(2,p)$ and let $H<G$ be a subgroup of index $p$. We are going to show that $p\le 11$.
We have $|G|=(p+1)p(p-1)/2$, hence $|H|=(p+1)(p-1)/2$.
The order of the stabilizer $H_\infty$ divides both the order of $H$ and the order of $G_\infty$.
The latter is the group of linear functions $\{x\mapsto ax+b\mid a,b\in\mathbb{F}_p, a\mathrm{\ is\ a\ square}\}$
of order $p(p-1)/2$, so $|H_\infty| \mid (p-1)/2$. On the other hand, the index of $H_\infty$
in $H$ is at most the number of permuted elements (``letters'' in Galois's terminology), that is
$|H:H_\infty|\le p+1$. Thus we obtain $|H_\infty|=(p-1)/2$ and $|H:H_\infty|=p+1$, hence
$H$ is transitive.
From the Lemma above we infer that
$H_\infty$ fixes some $c\in \mathbb{F}_p$. Let $t\in G$ be the translation $x\mapsto x+c$.
Then the conjugate subgroup $t^{-1}H_\infty t$ fixes $\infty$ and $0$. So, replacing $H$ by
$t^{-1}Ht$, we may assume without loss of
generality that $H_\infty$ fixes $0$ as well, so $H_\infty=H_0$.

Now we define an equivalence relation on $\{0,1,\dots,p-1,\infty\}$
by setting $x\sim y$ if $H_x=H_y$. This relation
is clearly $H$-invariant, since it follows that $H_{h(x)}=hH_xh^{-1}=hH_yh^{-1}=H_{h(y)}$, whenever
$x\sim y$ and $h\in H$, so $h(x)\sim h(y)$. We have seen that the equivalence class of $\infty$ is
$\{\infty,0\}$. If $x\sim y$ are different from $\infty$ and $0$, then for any square $a\in\mathbb{F}_p$
we can apply the transformation $x\mapsto ax$, which is now an element of the subgroup $H$,
and obtain $ax\sim ay$.
By the transitivity of $H$ it follows that every class of $\sim$ has two elements,
hence we will talk about partners of elements (linked letters --- ``lettres conjointes'' --- in Galois's manuscript).
The $H$-invariance of $\sim$ means that $H$ is an imprimitive permutation group.

Following Galois, let us denote by $M$ the partner of $1$. First we consider the case when
$M$ is a square. Then $1\sim M$ implies $M\cdot 1\sim M\cdot M$, hence $M^2$ is the partner of $M$,
so $M^2=1$, $M=-1$. (It can
happen only for primes $p\equiv 1\pmod 4$.) So $1\sim -1$, and then for every square $x$,
we have $x=x\cdot1\sim x\cdot(-1)=-x$. Now the partners of squares are squares themselves,
hence the partner of any non-square must be a non-square. If $x\sim M'x$ for non-squares $x$
and $M'x$, then $M'$ is a square, and it follows that
$M'=-1$ by the same reasoning as above. Summarizing, we
have obtained that if the partner of $1$ is a square, then the partner of every $x\notin\{\infty,0\}$ is
$-x$.

Now we look at the case when $M$ is not a square. Then for every nonzero square $x$ we have that
$x=x\cdot 1\sim x\cdot M$, so the partner of every nonzero square $x$ is the non-square $Mx$.

Now knowing the equivalence relation $\sim$ (depending on $M$, of course)
we have to check whether the subgroup of $G$ consisting
of those fractional linear transformations that preserve $\sim$ is indeed a transitive subgroup, as
$H$ has to be. So $H$ must contain a permutation $h$ that maps $1$ to $\infty$. Since $H$ preserves
$\sim$, we get that $h$ maps $M$ to $0$. Thus the fractional linear transformation $h$ must have
the form
\[
x\mapsto a\frac{x-M}{x-1}
\]
with an $a\in\mathbb{F}_p$ such that $a(M-1)$ is a square. Note that this $h$
must preserve the relation $\sim$.

In the first case, when $M=-1$ and the partner of every $x\notin\{\infty,0\}$ is $-x$, we must have for every
$x\notin\{\infty,0,1,-1\}$ that $h(-x)=-h(x)$, i.e.,
\[
a\frac{-x+1}{-x-1}=-a\frac{x+1}{x-1},
\]
hence $x^2+1=0$. This equation can hold for at most two elements $x$,
therefore we must have $|\{2,3,\dots,p-2\}|\le 2$, hence $p\le 5$ in this case.

In the second case, when the partner of every nonzero square $x$ is $Mx$ with a fixed non-square $M$, we
must have for every square $x\ne 0,1$ that $\{h(x),h(Mx)\}=\{y,My\}$ for some square $y$. It means
that either
\[
Ma\frac{x-M}{x-1}=a\frac{Mx-M}{Mx-1}
\]
or
\[
Ma\frac{Mx-M}{Mx-1}=a\frac{x-M}{x-1}.
\]
That is, for every square $x\ne 0,1$ either $x^2-(1+M)x+1=0$ or $x^2-\left(1+\frac{1}{M}\right)x+1=0$
holds. Each of these quadratic equations can have at most two roots in $\mathbb{F}_p$, so we obtain
that the number of non-zero squares in $\mathbb{F}_p$ is at most five, thus $p\le 11$.
\end{proof}

\section{Examples in the cases $p=5, 7, 11$}\label{examples}

We have seen that if a subgroup of index $p$ exists in $\PSL(2,p)$ then it is imprimitive,
namely,
it preserves a partition of $\{0,1,\dots,p-1,\infty\}$ into $(p+1)/2$ pairs. If we give an
appropriate partition then it is a matter of checking if the set of permutations in
$\PSL(2,p)$ that preserve the given partition is indeed a subgroup of index $p$ (or,
equivalently, of order $(p+1)(p-1)/2$). In fact, Galois gives the suitable partitions:

\hskip2cm for $p=5$: $\infty,0\mid1,4\mid2,3$ \ ($M=-1$);

\hskip2cm for $p=7$: $\infty,0\mid1,3\mid2,6\mid4,5$ \ ($M=3$);

\hskip2cm for $p=11$: $\infty,0\mid1,2\mid3,6\mid4,8\mid5,10\mid9,7$ \ ($M=2$).

(In the case $p=5$ it is not given explicitly, but since the pairs $\infty,0$ and $1,-1=4$ are mentioned,
the remaining two elements trivially yield the missing pair.)

In the last two cases he also writes down the elements of the subgroup in the form
\[
a \frac{x-b}{x-c},
\]
where $b,c$ is a pair, and $a$ has the same quadratic character as $c$. (Formally,
it should have the same quadratic character as $b-c$, as it was corrected in later
editions of Galois's works, but in the examples given for $p=7$ and $p=11$ it does
not matter. Furthermore, the meaning of $a\frac{x-\infty}{x-0}$ and $a\frac{x-0}{x-\infty}$
must be clarified. Namely, these should be $-a/x$ and $ax$, respectively, with $a$
being a square.)

Our discussion in Section~\ref{theproof} yields that if we choose $\infty, 0$ as one class of the partition, then
the other classes have either the form $x,-x$ (and this can happen only for $p=5$) or
$x,Mx$ for all nonzero squares $x$ with a fixed non-square $M$. In the latter case every square
$x\ne 0,1$ must satisfy one of the quadratic equations $x^2-(1+M)x+1=0$ or $x^2-(1+\frac{1}{M})x+1=0$,
furthermore, $p\le 11$. This necessary condition allows the values $M=3$ or $5$ for $p=7$, and
$M=2$ or $6$ for $p=11$. If some value of $M$ works then so does $1/M$ as well, because the
transformation $x\mapsto 1/x$ interchanges the two partitions, so if the subgroup
preserving the partition corresponding to $M$ has the required order, then conjugating
it by $x\mapsto 1/x$ yields a subgroup preserving the partition corresponding to $1/M$.
So in each case ($p=5,7,11$) there is essentially one suitable partition, namely, the one
found by Galois. (That is, all subgroups of index $p$ in $\PSL(2,p)$ are conjugate
in the larger group $\PGL(2,p)$.)

What has to be done is to show that the given fractional linear transformations
form a subgroup. Their number is indeed $(p+1)(p-1)/2$, since in $a(x-b)/(x-c)$ we can
choose $c$ arbitrarily (in $p+1$ ways), then it determines its partner $b$ uniquely,
and finally the quadratic character of $a$ should be the same as the one of $b-c$,
so we have $(p-1)/2$ choices for $a$. Galois probably verified it by calculating
compositions of such transformations. There are only very few traces of this in the
extant manuscripts in the holdings of the Library of the {\it Institut de France}
\cite{Galois}.

We give a concise verification using some geometric ideas. This approach is admittedly
anachronistic, since for Galois a group always meant a group of permutations of the
roots of an equation. The study of symmetry groups of geometric objects was a later
development in group theory (see \cite{Wussing}).

Our goal is to exhibit a subgroup $H_p$ of order $(p+1)(p-1)/2$ in $\PSL(2,p)$ for the
primes $p=5$, $7$, and $11$. So these subgroups should have order 12, 24, and 60, respectively.
These are the orders of the group of rotational symmetries of the Platonic solids.
In each case we will choose a $(p+1)$-element orbit of the symmetry group and label its
points by the elements of the projective line $\{0,1,\dots,p-1,\infty\}$ appropriately,
so that the permutations induced by the rotations on the chosen orbit will correspond to
elements of $\PSL(2,p)$. The orbit we take in the case of the regular tetrahedron will consist of
the six midpoints of the edges, for the cube and the icosahedron we take the set of
vertices. The suitable labelings are shown in the figures below.

\begin{center}
\includegraphics[height=47mm]{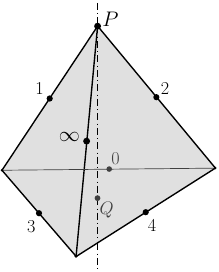} \qquad
\includegraphics[height=47mm]{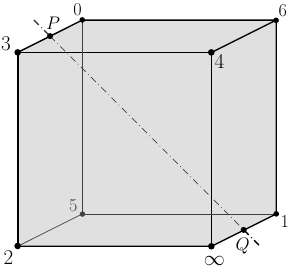} \qquad
\includegraphics[height=47mm]{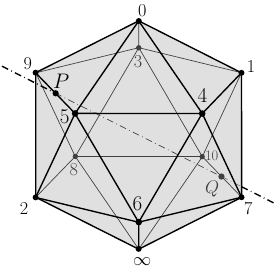}
\end{center}

To save calculations, it is enough to check that the permutations corresponding to
a pair of generators of the rotation group are indeed fractional linear transformations.
In the table below we give two rotations for each case. The axis of one of these rotations
goes through two (opposite) points (labelled by $0$ and $\infty$) of the orbit,
and the order of this rotation is $(p-1)/2$.
The other rotation does not fix any point in the chosen orbit, its axis is $PQ$, with
the points $P$ and $Q$ are marked in the figures. It is easy to see
that the two given rotations generate a transitive permutation group on the orbit of
size $p+1$, and since the stabilizer of a point contains a rotation of order $(p-1)/2$,
the subgroup generated by the two rotations is indeed the full rotation group of
order $(p+1)(p-1)/2$ of the Platonic solid.

\begin{center}
\begin{tabular}{|c|c|c|c|c|c|}
  \hhline{|=|=|=|=|=|=|}
  $p$ &polyhedron &axis &angle &permutation & function \\
  \hhline{|=|=|=|=|=|=|}
  5 &tetrahedron & $0\infty$ & $180^\circ$ & (0)(1 4)(2 3)($\infty$) &$4x$ \\
  \cline{3-6}
   & &$PQ$  &$120^\circ$ &(4 0 3)(1 $\infty$ 2) &$2(x-4)/(x-1)$ \\
  \hhline{|=|=|=|=|=|=|}
  7 &cube & $0\infty$ &$120^\circ$ & (0)(1 2 4)(3 6 5)($\infty$) & $2x$ \\
  \cline{3-6}
   & &$PQ$ &$180^\circ$ &(3 0)(1 $\infty$) (2 6) (4 5) &$(x-3)/(x-1)$  \\
  \hhline{|=|=|=|=|=|=|}
  11 &icosahedron & $0\infty$  &$72^\circ$ &(0)(1 3 9 5 4)(2 6 7 10 8)($\infty$) & $3x$ \\
  \cline{3-6}
  & &$PQ$ &$180^\circ$  &(2 0) (1 $\infty$) (3 6) (4 8) (5 9) (7 10) & $(x-2)/(x-1)$ \\
  \hhline{|=|=|=|=|=|=|}
\end{tabular}
\end{center}

A different approach was taken by John H.\ Conway \cite{Conway}, who considered Galois's partitions
as products of transpositions (($\infty$ 0) (1 4) (2 3) for $p=5$, etc.), and showed that each of them has
exactly $p$ conjugates under the action of $\PSL(2,p)$, so this set of involutions can be considered
as the required $p$ objects on which the the group $\PSL(2,p)$ ($p=5,7,11$) acts.

In another related work, Bertram Kostant \cite[p.~364]{Kostant} mentions that the set of vertices of the icosahedron
may be identified with the projective line over the 11-element field in such a fashion that there exists a partition
into six pairs of elements whose stabilizing subgroup of $\PSL(2,11)$ is $\PSL(2,5)$. Unfortunately, he does not
give the explicit identification in this paper, neither does he mention that this pairing can be found in
the last letter of Galois.

\medskip

{\bf Acknowledgements.} In October, 2020 I have exchanged a few e-mails with Peter M.\ Neumann
on this subject. His kind comments were helpful for my project. The nice figures of the three
polyhedra were produced by G\'eza K\'os.

\end{document}